\documentclass[]{bolmat5}
\usepackage{amsfonts}
\usepackage{amssymb}
\usepackage{amsmath}
\usepackage{tikz-cd}
\usepackage{float}
\usepackage{tikz}
\usetikzlibrary{trees}
\usetikzlibrary{decorations.pathmorphing}
\usetikzlibrary{decorations.pathreplacing}
\usetikzlibrary{decorations.markings}
\usetikzlibrary{decorations.shapes}
\usepackage{verbatim}

\usepackage{amssymb, eurosym}
\usepackage{amsmath}
\usepackage{textcomp}
\usepackage[includeheadfoot,a4paper]{geometry}
\usepackage{multicol}
\usepackage{amsfonts}
\usepackage{color}
\usepackage{enumerate}
\labeldocument[firstpage = 1, volume = 0, number = 0, month = 00, year = 1900, day = 00, monthreceived = 0, yearreceived = 1900, monthaccepted = 0, yearaccepted = 1900]

\begin{document}
	
\title[maintitle = {Geometric contexts and applications to logic},
	othertitle = {Contextos geométricos y aplicaciones a la lógica},
	shorttitle = {}
]
	
\begin{authors}[] % Campo opcional para colocar el encabezado de autores.
\author[firstname = {Mayk},
	surname = {Alves de Andrade},
	institutionnumber = {1},
	email = {andrade.mayk@ime.usp.br},
]
\author[firstname = {Hugo Luiz},
	surname = {Mariano},
	institutionnumber = {1},
	email = {hugomar@ime.usp.br},
]
\end{authors}
	
\begin{affiliations}
\affiliation[
	department = {Institute of Mathematics and Statistics},
	institution = {University of S\~ao Paulo},
	city = {S\~ao Paulo},
	country = {Brazil}
]
%\affiliation[
%	department = {Department or School},
%	institution = {University 2},
%	city = {City},
%	country = {Country}
%]
\end{affiliations}

\begin{mainabstract} 

There are many examples of dualities between topological spaces and algebras in the literature. Particularly, many of those examples come from the algebraic counterpart of a logical system, e.g, boolean and heyting algebras, MV-algebras, monoidal categories and linear logic etc. This process can be generalized to a representation of an algebraic structure by a sheaf in a suitable topological space. Recently, Awodey and his contributors constructed a duality between first order logic and topological grupoids from which they have proven some results of completeness. Furthermore, they built a notion of scheme of a first order theory. The attempt of this project is to use the construction of abstract schemes invented by Toen to give an unified approach from which we can generalize the aforementioned results.

\keywords{Sheaf Theory;  Geometric Contexts; Logic}
\end{mainabstract}
	
\begin{otherabstract} Hay muchos ejemplos de dualidades entre espacios topológicos y álgebras en la literatura. En particular, muchos de esos ejemplos provienen de la contraparte algebraica de un sistema lógico, por ejemplo, álgebras booleanas y heyting, MV-álgebras, categorías monoidales y lógica lineal, etc. Este proceso se puede generalizar a una representación de una estructura algebraica por un haz en un espacio topológico adecuado. Recientemente, Awodey y sus colaboradores construyeron una dualidad entre la lógica de primer orden y los grupoides topológicos a partir de la cual han probado algunos resultados de completitud. Además, construyeron una noción de esquema de una teoría de primer orden. El intento de este proyecto es utilizar la construcción de esquemas abstractos inventados por Toen para dar un enfoque unificado a partir del cual podamos generalizar los resultados antes mencionados.
		
\keywords{Teoria de haces; Contexto Geométricos, Lógica}
		
\end{otherabstract}
	
\msc{Mathematics Subject Classification According to AMS.}
	
%\section{First section}
\section{Dualities}

There are several examples of dualities between algebraic categories and topological/ geometrical spaces. The most famous ones are the Stone duality between boolean algebras and Stone spaces and the duality between commutative rings and affine spaces. We follow the discussion on \cite{Awodey}.

\subsubsection*{\textbf{Stone duality}:}

Let $B$ be a Boolean algebra. Consider the set $\textbf{Spec}(B)$ consisting of prime/maximal ideals of $B$, which are exactly the complements of ultrafilters. Endow $\textbf{Spec}(B)$ with a topology with a basis given by the open sets

\[B_f=\{p\in \textbf{Spec}(B) | f\notin p \}\]

Conversely, we can recover the Boolean algebra $B$ by considering the clopen spaces of $\textbf{Spec}(B)$. This process amounts to an equivalence of categories. Since clopen subspaces are continuous functions $f:\textbf{Spec}(B)\rightarrow \textbf{2}$, we get the \textit{Stone representation for boolean algebras}, i.e., there is always a injective homomorphism 

\[B \rightarrow \textbf{2}^{\textbf{Spec}(B)}\]

 from a Boolean algebra to the set of continuous functions on the local Boolean algebra $\textbf{2}$.

\subsubsection*{ \textbf{Grothendieck duality}:}

The Grothendieck duality is the milestone of algebraic geometry. Remember that a ring is called local if $1\neq 0$ and $x+y=1$ implies $x=1$ or $y=1$. 

\begin{theorem}
\textbf{Grothendieck ring representation} Let $A$ be a commutative ring. Then there is a space $\textbf{Spec}(A)$ and a sheaf of rings $\mathcal{R}$ such that:
\begin{enumerate}
    \item for every $p\in \textbf{Spec}(A)$, the stalk $\mathcal{R}_p$ is a local ring;
    \item there is an isomorphism $A\cong \Gamma (\mathcal{R})$, where $\Gamma(\mathcal{R})$ is the ring of global sections.
\end{enumerate}
The space $\textbf{Spec}(A)$ is built essentially in the same way as for Boolean spaces: the points are the prime ideals of $A$ and a base for a topology (the Zariski topology) is given by sets

\[A_f=\{p\in \textbf{Spec}(A) | f\notin p \}\]

The \textit{Grothendieck representation theorem}, therefore, says that for every commutative ring $A$ there is a injective homomorphism $A\rightarrow \Pi_p \mathcal{R}$ from $A$ to a product of local rings.
\end{theorem}
\section{Schemes}

Toen \cite{ToVa}, aiming to establish a comparison between analytic geometry and algebraic geometry (a categorical version of the famous GAGA theorem) defined a general procedure to define geometrical schemes. In the same way a manifold can be seen as a gluing of topological spaces through homeomorphisms, a geometrical scheme is a gluing of objects through a specified class of morphisms \textbf{P}. 
	\begin{definition}
	    %Um contexto geométrico é uma tripla $(\mathcal{C},\tau,\textbf{P})$ onde $\mathcal{C}$ é uma categoria, $\tau$ é uma topologia de Grothendieck em $\mathcal{C}$ e \textbf{P} é uma classe de morfismos de $\mathcal{C}$ de modo que:
     A \textbf{geometrical context} is a triple $(\mathcal{C},\tau, \textbf{P})$ where $(\mathcal{C},\tau)$ is a site and \textbf{P} is a class of morphisms such that 
	    \begin{enumerate}
	        \item the topology $\tau$ is \textbf{subcanonical};
	        \item \textbf{P} is \textbf{admissible}: it has all the identities, is stable under pullbacks and closed under composition;
	        \item \textbf{P} is $\tau$\textbf{-local}: if $\varphi:U\rightarrow V$ is a map and there is a covering $(\rho_i: U_i\rightarrow U)_{i\in I}$ such that $\rho_i\in \textbf{P}$ e $\varphi\circ\rho \in \textbf{P}$ for all $i\in I$, then $\varphi\in \textbf{P}$;
	        \item $\tau$ is \textbf{P}\textbf{-generated}: for all $U$ and all coverings $(\rho_i: U_i\rightarrow U)_{i\in I}$, there is a family of morphisms $(\varphi_j:V_j\rightarrow U)_{j\in J})$, with $\varphi_j\in \textbf{P}$ for all $j\in J$, such that each $\rho_i:U_i\rightarrow U$ factorizes through $\varphi_j:V_j\rightarrow U$ for some $j\in J$;
	        \begin{center}
	         % https://q.uiver.app/?q=WzAsMyxbMCwwLCJVX2kiXSxbMSwwLCJVIl0sWzAsMSwiVl9qIl0sWzAsMSwiXFxyaG9faSJdLFsyLDEsIlxcdmFycGhpX2oiLDJdLFswLDJdXQ==
\begin{tikzcd}
	{U_i} & U \\
	{V_j}
	\arrow["{\rho_i}", from=1-1, to=1-2]
	\arrow["{\varphi_j}"', from=2-1, to=1-2]
	\arrow[from=1-1, to=2-1]
\end{tikzcd}
	   \end{center}
    
        \item \textbf{P} is \textbf{locally cartesian}: if $\phi: U\rightarrow V$ is a \textbf{P}-morphism, then there is a covering $(\rho_i:U_i\rightarrow U)_{i\in I}$, with $\rho_i \in \textbf{P}$ for all $i\in I$, such that all morphisms $\phi\circ \rho_i$ are cartesian.
	        
	        \end{enumerate}
	\end{definition}
 The central construction in a geometrical context is that of a scheme. 
 \begin{definition}
     Let $(\mathcal{C},\tau, \textbf{P})$ be a geometrical context. A scheme is a sheaf $F$ in $(\mathcal{C},\tau)$ such that there exists a family of objects $\{U_i\}_{i\in I}$ and a morphism $p:\coprod_{i\in I}h_{U_i}$ satisfying the following conditions:
     \begin{enumerate}
         \item $p$ is a sheaf epimorphism;
         \item the morphism $h_{U_i}\rightarrow F$ is an open immersion (i.e., represented by a family of affine schemes $h_X$).
     \end{enumerate}
     The objects $U_i$ together with the map $p$ will be called an \textbf{open atlas} for $F$.
 \end{definition}
 In this way, we have a notion of scheme for any geometrical context. Furthermore, we can construct a notion of translation between geometrics contexts. Following \cite{Vaq2021}, let $(\mathcal{C},\tau,\textbf{P})$ and $(\mathcal{D},\rho,\textbf{Q})$ be two geometrical contexts, and let $f : C \rightarrow D$ be a functor. Assume that the functor $f$ satisfies the following
conditions:
\begin{enumerate}
    \item The functor $f$ commutes with finite limits and is continous for the topologies $\tau$ and $\rho$.
    \item The image by $f$ of a morphism in $P$ is a morphism in $Q$.
\end{enumerate}
Then the morphism of sites $\bar f_! : Sh (C) \rightarrow Sh(D)$
preserves schemes and sends morphisms in $P$ into morphisms in $Q$.

The analytification functor involved in the equivalence between analytical geometry and algebraic geometry, the GAGA theorem of Serre, is such a functor.

\section{Stone type dualities}
The list of Stone type dualities are extensive. We found in the literature some attempts to give general treatments of such constructions. Johnstone's book has a great historical account on this subject, as well as many references. We chose two references which are indicated there and one which follows directly from the ideas of the book.

\subsection{Universal algebra}

Davey \cite{Davey} gave a procedure to construct duality for varieties of algebras. This is the earliest paper we could find approaching the subject using the tools of universal algebra. 

A \textbf{congruence} of an algebra $A$ is an equivalence relation $\theta$ on its carrier $A$ such that for any operation $\omega$ of $A$:
\begin{center}
    if $\theta (a_1,b_1),\dots, \theta (a_n,b_n)$, then $\theta(\omega(a_1,\dots,a_n), \omega(b_1,\dots,b_n))$.    
\end{center}
\begin{theorem} 
  \textbf{\cite{Davey}}  Let $(\Theta_i)_{i\in I}$ be a family of congruences on an algebra $A$, and assume $A$ is a subdirect product of the algebras $(A/\Theta_i)_{i\in I}$. Then there exists a topological space $X$ and a sheaf $\mathcal{F}_A$ such that:

    \begin{enumerate}
        \item The stalk $\mathcal{F}_x$ is (isomorphic to) $A/\Theta_x$;
        \item the homomorphism $p:A \rightarrow \Gamma\mathcal{F}(X)$ embeds $A$ into the algebra of global sections of $\mathcal{F}_A$.
    \end{enumerate}
    In this case, the space $X$ is just the collection $\{\Theta_x\}$ of congruences endowed with the topology given by the subbasis $\{\Theta_i | (a,b) \in \Theta_i\}$ for any $a,b\in A$.
\end{theorem}

 Many of the sheaf representation of algebras found in the literature are instances of the above construction. Knoebels book gives a detailed account of it, as well as developments and applications. Sheaves over Boolean algebras, for instance, are well-known to be forcing models of set-theory.

A close but slightly more general definition is given by Caicedo (\cite{caicedo}), following Macintyre (\cite{mac}) and later explored by Villaveces and Padilla (\cite{villa}). They define a sheaf of \textbf{structures} in a first-order language.

\begin{definition}
    Consider a first-order language $L$. Then a sheaf of $L$-structures $\mathcal{F}_L$ is a topological space $X$ together with:
    \begin{enumerate}
        \item a sheaf $(E,p)$ over $X$ç
        \item for each $x\in X$, a $L$-structure $(E_x, R_x,\dots,f_x,\dots,c_x,\dots)$ for the fiber $E_x$. By abuse of notation, we call such structure simply by $E_x$.

        plus some continuity conditions:
        \begin{enumerate}
            \item the set $R=\cup_xR_x$ is an open of $\cup_x E_x^n$;
            \item the map $f:\cup_xf_x:\cup_x E_x^m\rightarrow \cup_x E_x $ is continuous
            \item the map $c:X\rightarrow E$ given by $c(x)=c_x$ is a continuous section.
        \end{enumerate}
    \end{enumerate}
\end{definition}

\subsubsection{Topos-theoretic approach}

The above definition closely resembles Butz \cite{Butz} construction. In a short historical perspective, this construction goes back to \cite{Grothendieck} where Grothendieck showed that the etale topos of a field $k$ is equivalent to the category of continuous $G$-sets for $G$ the profinite Galois group of a separable closure of $k$.

Joyal and Tierney in \cite{JT} then developed a much more general theorem, showing that any Grothendieck topos is equivalent to the category $Sh_G(X)$ of equivariant sheaves over a specific localic groupoid $G$ (a groupoid internal to the category of locales). As a lemma, they proved also that every topos $\mathcal{E}$ can be covered by a category of sheaves $Sh(X)$ over some locale $X$. Makkai and Barr in \cite{MakkaiBar} say that the covering result in \cite{JT} was proved already in \cite{Makkai} and draw the attention to the fact that "model-theoretic methods as such (the ones introduced in \cite{Makkai}) remained foreign to the
topos literature. Topos theorists tend to ignore results that are proved
using model theory."

To understand the preceding affirmation, recall that every Grothendieck topos $\mathcal{E}_{\mathbb{T}}$ is the classifying topos of a geometric theory $\mathbb{T}$, that is, a (infinitary) first-order theory where the sentences are built using only finite conjuction $(\land)$, infinitary disjuction $(\bigvee)$ and existential quantifiers $(\exists)$ (see \cite{Makkai}). Notice how this connectives are easily interpreted using topology. The classifying topos comes from a site $(C_{\mathbb{T}},J)$ constructed using only the syntax of (infinitary) first-order logic: the objects are formulas, morphisms are classes of formulas provable in $\mathbb{T}$ and the topology $J$ is given by conditions described by provability relations). Moreover, models of $\mathbb{T}$ in any Grothendieck other Grothendieck topos $\mathcal{F}$ are in bijection with geometric morphisms $\varphi: \mathcal{F}\rightarrow \mathcal{E}_{\mathbb{T}} $. 
\[\text{Hom}(\mathcal F,\mathcal{E}_{\mathbb{T}})\simeq \mathbb{T}-\text{mod}(\mathcal F)\]

where a geometric morphism $\varphi$ is given by an adjunction $\varphi^*\vdash \varphi_*$ where the left adjoint $\varphi^*$ preserves finite limits.

Particularly, the category of sets $\textbf{Set}$ is a Grothendieck topos and, therefore, set-theorical models of a geometric theories are also classified by $\mathcal{E}_{\mathbb{T}}$. In addition, the etale topos is the classifying topos for the theory of separably closed extensions of $k$.

In \cite{herman} it was proved that any Turing machine can be encoded in the category of functors from $\omega$ to finite sets. This category forms a pretopos of calculations.

Finally, Butz in his PhD thesis took the assertion made by Makkai and Reyes seriously and developed a representation theorem for toposes with enough points using for this a description which relies deeply in a model theoretical framework. Explicitly, he proves the following theorem: 

\begin{theorem}
    (\cite{Butz2} and \cite{Butz3}) Let $\mathcal{E}$ be a topos with enough points. There exists a topological groupoid $G\rightrightarrows X_\mathcal{E}$ and a topos morphism $\varphi: Sh(X_{\mathcal{E}} ) \rightarrow \mathcal{E}$ , such that
    \begin{enumerate}
        \item $\varphi$ is a full and faithful embedding of $\mathcal{E}$ into $Sh(X_{\mathcal{E}})$ which induces an equivalence of topoi \[\mathcal{E}\simeq Sh_G(X).\] ;

        \item furthermore, for any Abelian group $A$ in $\mathcal{E}$, the morphism $\varphi$ induces isomorphisms between the cohomology groups \[H^n(\mathcal{E},A)\simeq H^n(X_{\mathcal{E}},\varphi^*(A))\]
\end{enumerate}
\end{theorem}

 They also are valid for the localic. Adapting this construction and inpired by Moerdjik, Forssell in \cite{forssell} proves that it can be extend to an adjunction btween theories and contiuous groupoids over the empty theory groupoid:

% https://q.uiver.app/#q=WzAsMixbMCwwLCJcXG1hdGhjYWx7VH0iXSxbMiwwLCJcXHRleHRiZntHcnB9L1xcbWF0aGNhbHtHX1xcdmFybm90aGluZ30iXSxbMSwwLCJGIiwyLHsiY3VydmUiOjJ9XSxbMCwxLCJNIiwyLHsiY3VydmUiOjJ9XSxbMiwzLCIiLDIseyJsZXZlbCI6MSwic3R5bGUiOnsibmFtZSI6ImFkanVuY3Rpb24ifX1dXQ==
\[\begin{tikzcd}
	{\mathcal{T}^{\text{op}}} && {\textbf{Grp}/\mathcal{G_\varnothing}}
	\arrow[""{name=0, anchor=center, inner sep=0}, "F"', bend right=30, from=1-3, to=1-1]
	\arrow[""{name=1, anchor=center, inner sep=0}, "M"', bend right=30, from=1-1, to=1-3]
	\arrow["\dashv"{anchor=center, rotate=-90}, draw=none, from=0, to=1]
\end{tikzcd}\]

It restricts to an equivalence for the case where the theories have enough $\textbf{Set}$ models and their semantical groupoids. This adjunction also restricts to the special case where the theory is \textbf{coherent}, that is, it is built using finitary formulas. This result is nothing but a \textbf{Stone type duality for first-order logic}!

A proposed generalization is given by Caramello \cite{Caramello}, explicitly following Johnstone's book, using the lattice of subobjects of an arbitrary (elementary) topos to give an unified construction of sheaf representations.
\begin{definition}
    Let $\mathcal{E}$ be a locally small cocomplete topos, $\Gamma$ be a subframe of $\textbf{Sub}_\mathcal{E}(1)$ and $i : X \rightarrow P$ be an indexing function of a set $P$ of points of $\mathcal{E}$ by a set $X$. The $\Gamma$-subterminal topology $\tau^\mathcal{E}_{\Gamma,i}$ on the set $X$ is the image of the function $\phi_{\Gamma,\mathcal{E}} : \Gamma \rightarrow \mathcal{P}(X)$ given by
    
\[\phi_{\Gamma,\mathcal{E}}(u) = \{x \in X |\ \varepsilon(x)^ *(u) \cong 1_{\textbf{Set}}\}\]

\end{definition}
\begin{proposition}
Let $\mathcal{C}$ be a preorder and $J$ be a Grothendieck topology on it. Then the topological space $X_{\tau,\textbf{Sh}(\mathcal{C},J)}$ is homeomorphic to the space which has as set of points the collection $\textbf{F}_{\mathcal{C}}^J$ of the $J$-prime filters on $\mathcal{J}$ and as open sets the sets the form \[\textbf{F}_I = \{ F \in \textbf{F}_\mathcal{C}^J / F \cap I  = \varnothing \},\]
where $I$ ranges among the $J$-ideals on $C$. In particular, a sub-basis for this
topology is given by the sets $\textbf{F}_c = \{ F \in \textbf{F}_\mathcal{C}^J / c\in F \}$, where $c$ varies among the elements of $C$.
\end{proposition}
\begin{proposition}
Let $\mathcal{C}$ be a preorder and $J$ be a Grothendieck topology on $\mathcal{C}$. Then the toposes $\textbf{Sh}(\mathcal{C},J)$ and $\textbf{Sh}(Id,J(\mathcal{C}))$ are equivalent.
\end{proposition}

Our first line of research to be proposed here is the investigation of a general and unified framework for model theory and related subjects. There are some proposals in the literature, such as institutions, accessible categories and ionads. In the same prospect, we want to establish what kind of logical systems can be given a Stone type representation. General notions of logic and translations between them have been studied in \cite{SusBr}, \cite{MenMa} and \cite{Pinto}. The first one is quite general and focuses in closure operators, which are just instances of subobjects of free sup-lattices. The second one focuses in the Blok-Pigozzi algebrization. The third one gives some first steps in the direction of the algebraic representation of logics, we want to look to it through the eyes of topology.

\subsection{Further connections}

In this section, we briefly comment another suggestive connections between model theory/logic and geometry. This connection have been sharpened since the pioneering work of Lawvere in his PhD thesis and subsequent publications. In \cite{Mac}, Macintyre comments how model since its beginning with Tarski work almost does not rely on set theory, but rather in the more general interplay between syntax and semantics. In fact, algebraic geometers and model theorists were working with similar problems but in different fashions.

For instance, the work of classification of Shelah gave rise to stability theory and then to geometric stability theory, or o-minimal geometry. Much of algebraic geometry has been translated to this setting, as we can see for instance in \cite{omin}, where a GAGA-type theorem has been demonstrated.

Another example is \cite{lascar}, where they show how the Lascar group of a first-order theory is isomorphic to the fundamental group of the space of models of a theory built out of models and elementary morphisms between them. In \cite{lascar2}, they show how to construct an abstract elementary class for each homotopy type.

In a similar manner, Arndt \cite{arnt} has recently applied abstract homotopy theoretical techniques to the abstract logic setting developed. More specifically, a model category structure for categories of logics and translations between them. Finally, all the homotopical interpretation of type theory which recently has led to a whole program to give new foundation to mathematics is nothing but another instance of the duality between syntax and semantics.

Finally, it is indispensable to mention Oostra's representation of Peirce's existential graphs in the complex plane , which suggests a logical formulation for each Riemann surface \cite{oostra}.

\section{Logical schemes}

As a final construction related to the interplay between algebraic geometry and first-order logic, there is the logical schemes constructed by Breiner following the earlier results we presented here.

A pretopos is a category which is both exact and extensive. They arise as the completion of coherent categories by the addition of coproducts and quotients. A pretopos with a power object is a topos.  Most importantly, every pretopos is equivalent to a coherent and equivariant sheaf on a topological groupoid. Such groupoid is built in the same wat we did before.
\begin{theorem}
\textbf{ \cite{Awodey}} Let $\mathcal{P}$ be a pretopos. There is a
topological groupoid $G\rightrightarrows X_\mathcal{P}$ with an equivariant sheaf of pretoposes $\widetilde{\mathcal{P}}$ such that: 
\begin{enumerate}
    \item for every $g \in G$, the stalk  $\widetilde{\mathcal{B}}_g$ is a well-pointed pretopos;
    \item for the global sections of $\widetilde{\mathcal{P}}$ there is an equivalence $P \cong \Gamma(\widetilde{\mathcal{P}})$.
Thus every pretopos is equivalent to the global sections of a sheaf of
local pretoposes.
\end{enumerate}
Consequently, for every pretopos $\mathcal{E}$, there is a pretopos embedding 
\[\mathcal{E}\rightarrow \Pi_g \mathcal{E}_g\]
with each $\mathcal{E}_g$ a local pretopos.
\end{theorem}
For more on the subject, see \cite{breiner} and \cite{rios}
\begin{definition}

    A \textbf{logical scheme} $(G\rightrightarrows X,\mathcal{O}_X)$ is a topological groupoid together with an equivariant sheaf of pretoposes such that every stalk is a local pretopos and there is an open cover $\{U_i\}$ by affine subspaces $U_i \cong X_{\mathcal{E}_i}$ for pretoposes $\mathcal{E}_i$.

    A morphism of logical schemes $(f,\varphi): (G\rightrightarrows X,\mathcal{O}_X)\rightarrow (H\rightrightarrows Y,\mathcal{O}_Y)$ is a morphism $f:X\rightarrow Y$ between the underlying groupoids and a pretopoi morphism $\varphi:\mathcal{O}_Y\rightarrow f_*\mathcal{O}_X$ internal to the cateory of coherent equivariant sheaves over $H\rightrightarrows Y$ such that $f^*$ is conservative on stalks.

    \end{definition}

    Naturally, both the representation and the logical schemes specialize to the Boolean case. Notice how this construction opens many possibilities to explore the relation between algebraic geometry and logic. For any theorem in the theory of schemes there is a corresponding analogue here. Geometrical contexts could offer a bridge between both worlds.

    \section{Algebraic geometry = geometrical tools for logic}

    Hodges in \cite{hodges} defines model theory as "algebraic geometry minus fields". We wanted with this review at this motto has had some serious unfoldings in past fifty years or so. Algebraic geometry and model theory share unavoidable links. The duality algebra/geometry is the same as the duality syntax/semantics and there are many situations to be explored. The whole tradition of anlytical philosophy remains to be fused with modern geometry. The result of this fusion will not be model theory or geometry, but an elevated version of both.

    As a final remark, we would like to note that a powerful tool which has not been well explored yet is the notion of cohesiveness, which played an important role through the works of Lawvere (e.g., in \cite{law1} and \cite{law2}). For instance, Schreiber (\cite{sche}) defined synthetically in any cohesive $\infty$-topos important notions of differential geometry. Of particular interest, it is possible to define synthetically what would be a connection on generalized principal bundles.

    Remember that principal $G$-bundles over a topological space $X$ are classified by geometrical morphism $Sh(X)\rightarrow BG$, where $BG$ is the topos of $G$-sets. But $G$-sets are nothing but an equivariant sheaf over the point. Therefore, by what we have seen previously, models of a geometric theory behave much like principal bundles for the topological groupoid corresponding to that theory. 
    
    The question we want to investigate, therefore, is the following: is there a specific kind of theories for which their classifying toposes and respective models give rise to a cohesive setting? In the analogy with usual geometry, what kind of information a connection would carry? What can we recover from modern formalism of physical theories presented at \cite{sche} in this setting?

    The aplications could be inumerous. For instance, a notion arising from symplectic geometry could be usefull to formulate computational constraints based on analogies with energy conservation in physical systems. Also, could share light in a new notion of automated reasoning. Even though this is only speculative, it is worthy to explore the full length of possibilities brought by the relations established so far.

\section{Conclusion}

    Based upon the developments aforementioned in this review, we hope that the reader is convinced about the viability of the use of tools of algebraic and differential geometry to reformulate problems in logic and related subjects. The usual cohomology and homotopy theories arising in geometry certainly are a source of invariants for logical system that could be related to provability, computability, definibility, normalization and so on. It is only a matter of time before new and deep advances take place in this branch of research.

     % % % % % % % % % % % % % % % % % % % % % % % % % % % % %
%\bibliography{BMBibTeX}

\begin{thebibliography}{10}

{\footnotesize

	\bibitem[{\bf ToVa08}] {ToVa} Toen, B., Vaquié, M., {\em Algébrisation des variétés analytiques complexes et catégories dérivées}, Math. Ann. 342, 2008, 789–831.

    \bibitem[{\bf Awo21}]{Awodey} Awodey, S. {\em Sheaf Representations and Duality in Logic}. In: Casadio, C., Scott, P.J. (eds) Joachim Lambek: The Interplay of Mathematics, Logic, and Linguistics. Outstanding Contributions to Logic, vol 20. Springer, Cham. 2021.

     \bibitem[{\bf Car11}] {Caramello} Caramello, O., {\em A topos-theoretic approach to Stone-type dualities}, arXiv:1103.3493, 2011.

     \bibitem[{\bf Dav73}]{Davey} Davey, B.A. {\em Sheaf spaces and sheaves of universal algebras}. Math Z 134, 275–290, 1973.

         \bibitem[{\bf But96}]{Butz} Butz, C. {\em Logical and cohomological aspects of the space of points of a topos}. PhD Thesis, 1996.

           \bibitem[{\bf BuMo98}]{Butz2} Butz, C. Moerdjik, I. {\em Representing topoi by topological groupoids}. Journal of pure and applied algebra, 1998.

                 \bibitem[{\bf BuMo99}]{Butz3} Butz, C. Moerdjik, I. {\em Topological Representation of Sheaf Cohomology of Sites}. Compositio Mathematica, 1999.

     \bibitem[{\bf AdaRo94}] {AdaRo} Adamek, J., Rosicky, J. Locally Presentable and Accessible Categories (London Mathematical Society Lecture Note Series). Cambridge: Cambridge University Press, 1994.

      \bibitem[{\bf Vaq2021}] {Vaq2021} Vaquie, M. {\em Sheaves and functor of points.} In: J. Catren, G. Anel, M. (eds.) New spaces in mathematics: formal and conceptual reflections. Cambridge: Cambridge University Press, 2021.

      \bibitem[{\bf SGA1}]{Grothendieck} Grothendieck, A. {\em Revetements etales et groupe fondamental SGA 1}, 
Springer Lecture Notes in Mathematics 224 (1971). 

\bibitem[{\bf JoTi84}]{JT}  Joyal A. and Tierney. M. {\em An extension of the Galois theory of Grothendieck.}
American Mathematical Society, 1984.

 \bibitem[{\bf MaRe77}]{Makkai} Makkai, M. and Reyes. G. {\em First order categorical logic: Model-theoretical
methods in the theory of topoi and related categories.} Springer-Verlag, 1977.

\bibitem[{\bf Cai95}]{caicedo}  Caicedo, X.{\em Logica de los haces de estructuras.} 1995.


\bibitem[{\bf PaVi16}]{villa}  Padilla, G. and Villaveces, A. Non standard cohomology for equivariant sheaves: the role of generic models. 2016.


\bibitem[{\bf BaMa87}]{MakkaiBar} {\em On representations of Grothendieck toposes}
Canadian Journal of Mathematics, 1987.

\bibitem[{\bf For08}]{forssell}  {\em Henrik Forssell. First-Order Logical Duality}. PhD thesis, Carnegie Mellon
University, 2008.

\bibitem[{\bf SusBr73}]{SusBr} Susko, R. and Brown, D. {\em Abstract logcis and Classical abstract logics }. Dissertationes Mathematicae: Warsaw, 1973.

\bibitem[{\bf MenMa16}]{MenMa} Mendes, C. and Mariano, H. {\em Towards a good notion of categories of logics}. 2016.

\bibitem[{\bf Pin16}]{Pinto} Pinto, D. {\em A categorical foundation for the representation theory of logics}. 2016.

\bibitem[{\bf Mac73}]{mac}  Macintyre, A. Model completeness for sheaves of structures. Fundamenta Mathematicae, 1973.
\bibitem[{\bf Mac03}]{Mac}  Macintyre, A. {\em Geometrical and Set-Theoretic Aspects and Prospects.} The Bulletin of Symbolic Logic, 2003.

\bibitem[{\bf CCY18}]{omin} Campion, T., Cousins, G., Ye, J. {\em Classifying spaces and the Lascar group.}  2018.

\bibitem[{\bf CaYe19}]{lascar}  Campion, T., Ye, J. {\em Homotopy Types of Abstract Elementary Classes.}  2019.

\bibitem[{\bf BBT18}]{lascar2}  Bakker, B., Brunebarbe, Y., Tsimerman, J. {\em o-minimal GAGA and a conjecture of Griffiths.}  2018.

\bibitem[{\bf Arnt}]{arnt}  Arnt, P. {\em Private conversations }.

\bibitem[{\bf Oos19}]{oostra}  Oostra, P. {\em Complex Representation of Alpha Graphs for Implicative Logic
with Conjunction}. Boletín de Matemáticas, 2019.

\bibitem[{\bf Law06}]{law1}  Lawvere, W. {\em Some thoughts on the future of category theory
}. 2006.

\bibitem[{\bf Law07}]{law1}  Lawvere, W. {\em Axiomatic cohesion}. Theory and Applications of Categories, 2007.



\bibitem[{\bf Sch13}]{sche}  Schreiber, U. {\em Differential cohomology in a cohesive infinity-topos
}. 2013.

 \bibitem[{\bf Br14}]{breiner}  Breiner, S. {\em Scheme representation for first-order logic}. PhD Thesis, 2014.

 \bibitem[{\bf RiMa23}]{rios}  Rios, G., Mariano, H. {\em Model Theory Inspired by Grothendieckian Algebraic
Geometry: a Survey of Sheaf Representations for Categorical
Model Theory}. Latin American Journal of Mahtematics, 2023.

 \bibitem[{\bf Ho97}]{hodges} Hodges, H. A Shorter Model Theory. Cambridge University Press,
1997.

 \bibitem[{\bf Joh82}]{johnstone} Johnstone, P. Stone Spaces. Cambridge University Press,
1982.

\bibitem[{\bf MHG16}]{herman} Macedo, H., Haeusler, E. and Garcia, A. {\em Defining effectivenes unisg finite sets: a study on computability}. Cambridge University Press,
1982.




    %bibitem[Bel]{Bel} John L. Bell. \textit{Set theory: Boolean-valued models and independence proofs}. 3rd ed. Oxford Logic Guides, vol. 47. Oxford, United Kingdom: Clarendon Press, 2005.
    
    %\bibitem[Bel88]{Bel88} John L. Bell. \textit{Toposes  and  local  set  theories:  an  introduction}. Oxford  Logic Guides, vol. 14. Oxford, United Kingdom: Clarendon Press, 1988.
    
   % \bibitem[Bor08a]{Bor08a}  Francis  Borceux. \textit{Handbook  of  Categorical  Algebra}.  Vol.  1. Encyclopedia of mathematics and its applications, vol. 50. Cambridge, United Kingdom: Cambridge University Press, 2008.
    
  %  \bibitem[Bor]{Bor}  Francis  Borceux. \textit{Handbook  of  Categorical  Algebra}.  Vol.  3. Encyclopedia of mathematics and its applications, vol. 52. Cambridge, United Kingdom: Cambridge University Press, 2008.

   % \bibitem[BP]{BP} W. J. Blok, D. Pigozzi, {\em Algebraizable logics},  Memoirs of the  AMS {\bf 396}, American Mathematical Society, Providence, USA, 1989.


%\bibitem[Diac]{Diac} R. Diaconescu. {\em Institution-independent Model Theory}, { Birkhauser Basel · Boston · Berlin}, 2000.
%An Introduction to Partially Ordered Structures and Sheaves
%Francisco Miraglia
%Polimetricas




}

\end{thebibliography}
% % % % % % % % % % % % % % % % % % % % % % % % % % % % %

\end{document}